\input amstex
\documentstyle{amsppt}
\document
\magnification=1200
\NoBlackBoxes
\nologo
\pageheight{18cm}


\smallskip


\bigskip
 
\centerline{\bf MODULI, MOTIVES, MIRRORS\footnotemark1}
\footnotetext{Plenary talk at the 3rd European Congress
of Mathematicians, Barcelona, July 10--14, 2000}

\bigskip

\centerline{\bf Yuri I. Manin}

\medskip

\centerline{\it Max--Planck--Institut f\"ur Mathematik, Bonn}

\bigskip

\centerline{\bf 0. Introduction}

\medskip

This talk is dedicated to various aspects of Mirror Symmetry.
It summarizes some of the developments that took place since
M.~Kontsevich's
report [Ko2] at the Z\"urich ICM and provides an extensive,
although not complete bibliography.

\medskip

{\bf 0.1. Brief history.} Mathematical history of
Mirror Symmetry started in 1991, when an identity of a new type
was discovered in the ground--breaking paper by four physicists [CaOGP]
(it was reproduced in [MirS1] where  earlier works are
also described and motivated).

\smallskip

The left hand side (or $A$--side) of this identity was a generating series
for the numbers $n(d)$ of rational curves of various degrees $d$
lying on a smooth quintic hypersurface in ${\bold P}^4$.
The right hand side ($B$--side) was a certain hypergeometric function.
The Mirror Identity states that the two functions become identical
after an explicit change of variables which is defined
as a quotient of two hypergeometric functions of the same type. 

\smallskip

\smallskip

At the moment of discovery, not only the identity itself
remained unproved, but even its $A$--side was not well defined:
the correct way of counting rational curves was proposed
by M.~Kontsevich ([Ko4]) only in 1994. In the same remarkable paper
Kontsevich gave an explicit formula for $n(d)$ creatively using
Bott's fixed point formula for torus actions at the target space.
After the appearance of 
this paper
one could hope that the Mirror Identity for quintics (and
more general toric submanifolds)
ought to be provable by algebraic manipulations with
both sides. This turned out to be a difficult problem.
A.~Givental brought this program to a successful completion
in 1996, by introducing a new torus action at the
source space, stressing equivariant cohomology
and inventing ingenious calculational strategy (see [Giv2],
[Giv5], [Giv7], [BiCPP], [Pa]). For subsequent
important developments, see [LiLY1], [LiLY2], [Ber].

\smallskip

This work however did not unveil the mystery of
the Mirror Identity. The point is that the identity itself
was discovered by the physicists as only one
manifestation of a deeper principle. 
Physicists believe that  with any Calabi--Yau manifold $X$ one can associate
two $N=(2,2)$ Superconformal Field Theories (SCFT) which are the respective
$A$ and $B$ models (see e.g. [Wi1]). The Mirror Correspondence
between $X$ and $Y$ supposedly interchanges their $A$ and $B$ models.
In particular, in the case of quintics the hypergeometric
functions involved are actually periods of the mirror partner family
of our quintics, and $B$--models generally reflect properties
of variations of periods and Hodge structures.

\smallskip

Unfortunately, a precise and complete mathematical
definition of what constitutes an $N=(2,2)$ ScFT
is still lacking. Various components of this structure
with varying degree of precision are described
in the papers collected in [MirS1] and [MirS2]. 
In particular, a part of this structure is a modular
functor in the sense of Segal, with possibly infinite dimensional
Hilbert space. In turn, such theories are often constructed
via representation theory of a vertex algebra. See [MalSV] and 
[Bor2], [Bor3]  for the most recent mathematical approach to this picture,
achieving at least the construction of what seems to be the right
vertex algebra.

\smallskip

The parts that are involved in the statement of Mirror Identity above
refer correspondingly to the Quantum Cohomology
($A$--model, physicists' $\sigma$--model) and extended
variations of Hodge structure. Both are now well understood
mathematically: see [Ma5] and [Bar] respectively.
However, the Mirror partners are  connected
by much more ties than a mere Mirror Identity.
These ties, in particular, relate Lagrangian and complex geometry
in a remarkable way: see [StYZ] and [Ko2] for the basic conjectures
to this effect.

\smallskip

Therefore now, more than decade after it was discovered, the Mirror Symmetry
mathematically looks like a complex puzzle, some of the pieces
of which have found their respective places, some are still lying in disorder,
and some, most probably, are missing.

\medskip

{\bf 0.2. Plan of the paper.} This puzzle metaphor guided the organization of this report.

\smallskip

Section 1 is devoted to the binary relation
of {\it mirror partnership} between  
families of Calabi--Yau manifolds endowed with additional structures
which we call here
{\it cusps}. This relation consists in the
isomorphism of two Frobenius manifolds, constructed by
two different ways for the respective families. In turn,
Frobenius manifold isomorphisms generalize the Mirror Identity of
[CaOGP].

\smallskip

Section 2 explains various versions of another
mirror partnership relation, this time between certain 
symplectic, on the one hand, and
complex, on the other hand, manifolds, endowed 
with additional structure which in this case
is a choice of {\it a fibration by real tori}. Here I have
took
as starting point a part of Kontsevich's package [Ko2],
with further detalization taken from [StYZ], [PoZ], [AP],
and other papers. I have chosen for these relations the word
``partnership'', or ``duality'', as opposed to ``symmetry'', because 
the definition of both of them
is explicitly un--symmetric.

\smallskip

The next part of the section 2 restores the idea of 
Mirror Symmetry: according to [StYZ], the Mirror Symmetry
relation connects Calabi--Yau manifolds endowed with K\"ahler
structure, and thus simultaneously with compatible Lagrangian
and complex structures, so that Kontsevich's duality
can be imposed simultaneously upon two crossover Lagrangian/complex
pairs.  

\smallskip

Although the Frobenius manifold duality and the Lagrangian/complex
duality some day are expected to become parts of
a unified picture, at present contours of the latter
are rather vague.

\smallskip

One common part of both dualities is the prediction
of {\it mirror isomorphisms} connecting the cohomology
spaces of mirror partners $(X,Y)$. In particular, isomorphisms
of the Frobenius manifolds restricted to their spaces
of flat vector fields
produces isomorphisms $\mu_{X,Y}:\,
H^*(X,\bold{C})\to H^*(Y,\wedge^*(\Cal{T}_Y)\otimes V_Y^{-2}$
where $V_Y=H^0(Y,\Omega_Y^{\roman{max}}).$ 

\smallskip

Actually, algebraic geometric model of the $A$--side,
the theory of Gromov--Witten invariants, endows
$H^*(X,\bold{C})$ with much stronger structure,
which is motivic by its nature. It would be very
interesting to understand the geometry of the 
mirror reflection of Calabi--Yau motives.

\smallskip

However, even when the intricate inner workings of Mirror
Symmetry are understood, this will not be
the end of the story. 

\smallskip

{\bf 0.3. Dualities in string theory.} All this machinery emerged 
as an approximation in the quantum
superstring theory whose aim is to provide a unified
theory of matter and gravity (space--time). The first superstring revolution (1984--85) led to the belief that there are 
five consistent (perturbative, without ultraviolet
divergencies) superstring theories, each on ten--dimensional
space--time, or in other words,
10d Poincar\'e--invariant vacuum. For all of them, the low--energy
approximation is an effective 10d supergravity theory.
The second superstring revolution (1994--??) started with the Witten's suggestion that
all five theories are limits of a single theory (see review in
[Schwa]). In other words, they are perturbarive expansions
of a single underlying theory about distinct points
in the moduli space of quantum vacua. Moreover, a sixth special point
in this space is a 11d Poincar\'e--invariant vacuum.

\smallskip

C.~Vafa suggests to look at the
underlying {\it M--theory} as patched up from five/six
local descriptions and their compactifications, 
im much the same way as a manifold is patched up
from coordinate neighborhoods.  The transition
functions are called {\it dualities}, and according to
[StYZ], Mirror Duality is one of them.

\smallskip

If this is true, the Mirror Symmetry acquires an incredibly
high epistemological status as one of the building blocks
of the ambitious Unified Quantum Superstring Theory.

\smallskip

For mathematicians, this means that the puzzle we are trying to assemble,
is only a small piece of the still larger puzzle
whose contours are yet barely visible.

\newpage

\centerline{\bf 1. Frobenius manifolds and mirror partnership}

\smallskip

\centerline{\bf  between families of Calabi--Yau manifolds}

\medskip

{\bf 1.1. Calabi--Yau manifolds.} In this section I will
call a Calabi--Yau (CY) manifold {\it in the weak sense} any
projective (or compact K\"ahler) complex manifold $X$ with trivial canonical
sheaf. Any such manifold admits a finite unramified covering
$\widetilde{X}$ with the following property:

\smallskip

{\it $\widetilde{X}$ is a direct product of a complex torus
$A$, of a simply connected CY $Y$ with $h^{2,0}(Y)=h^{2,0}(\widetilde{X})$,
and of a simply connected CY $Z$ with $h^{2,0}(Z)=0.$}

\smallskip

If the factors $A$ and $Y$ are absent in any finite unramified
covering of $X$, then $X$ is CY {\it in the strong sense}.

\smallskip

In dimension 1, the only CY manifolds are elliptic curves,
in dimension 2, besides complex tori, there are
$K3$--surfaces. In dimension 3, the first examples of
CY in the strong sense appear. Quintics in $\bold{P}^4$
are the simplest of them. 
\smallskip

More generally, anticanonical
hypersurfaces in any compact toric manifold associated
with a {\it reflexive polyhedron} are CY's: see [Ba1].
This method produces 4319 families of $K3$--surfaces
and 473 800 776 families of CY threefolds,
among which at least 30 178 families can be distinguished
by their Hodge numbers (see [KrS]). It is still unknown,
whether the number of maximal families of CYs
in any dimension $\ge 3$ is finite or not.

\smallskip
The toric construction
(and its generalization to complete intersections in arbitrary
Fano manifolds)
remain the most important testing ground for basic
conjectures about CYs.

\smallskip

A general approach to the complex moduli spaces
of CY manifolds is furnished by the deformation theory.
The Kodaira--Spencer local versal deformation of
an $n$--dimensional CY manifold $X$ is unobstructed
and has dimension $h^{1,n-1}(X).$

\medskip

{\bf 1.2. Mirror partner families of CYs: preliminarities.}
The notion which we will describe in this section
is interesting mainly for CYs in the strong sense.
It develops the discovery made in [CaOGP].

\smallskip

This notion is an asymmetric  binary relation between versal local families 
$\{X_s\,|\,s\in S\}$,  $\{Y_t\,|\,t\in T\}$ of
CYs, satisfying the condition
$h^{1,1}(X)=h^{1,n-1}(Y)=r$ and endowed with some additional structure.

\smallskip

On the $A$--side, this additional structure consists
in a choice of a basis 
$(\beta_1,\dots ,\beta_r)$
of the group of numerically effective classes in $A_1(X_s)$.
When $s\in S$ varies,
elements of this basis must be horizontal with respect
to the Gauss--Manin connection. Such a basis determines
functions $q_j^A$ on $H^2(X_s,\bold{C}/\bold{Z})$: $q_j^A(L):=e^{2\pi i(L,\beta_j)}$.
We will refer to  $H^2(X_s,\bold{C}/\bold{Z})$ together with $q_j^A$
as {\it a K\"ahler cusp}.  
\smallskip

On the $B$--side, this additional structure consists
in the choice of a partial compactification $T\subset \overline{T}$
looking locally like the embedding of a  product
of pointed open unit discs in $\bold{C}$ into the product
of non-pointed unit discs. The variation
of Hodge structures of the family $Y_t$ must have
maximal unipotent monodromy on this compactification:
see [Mo2]. [Mo3], [De3]. Geometrically, $\overline{T}$
contains a point of ``maximal degeneration'' of the family
$Y_t$ and $T$ parametrizes Calabi--Yau manifolds ``with
large complex structure''.  The point
of maximal degeneration is the transversal intersection
of discriminantal divisors. Building upon [CaOGP] and [Mo1],
[Mo2], Deligne has shown in [De3] how to define a system
of functions $q_j^B$ on $\overline{T}$ in terms of
the variation of Hodge structures determined by $Y_t$. 
We will refer to the germ
of $\overline{T}$ at its point of maximal degeneration
as {\it moduli cusp} of the relevant moduli space.

\smallskip

The relation of {\it mirror partnership}
between such enhanced families $X/S$ and $Y/T$, in particular, 
identifies $q_j^A$ with $q_j^B$ and thus
establishes an isomorphism between a domain
in $H^2(X,\bold{C}/\bold{Z})$ where $q_j^A$ are sufficiently small
and the respective domain in the moduli cusp.
This isomorphism must identify two functions:
{\it potential of the small quantum cohomology}
at the $A$--side, and an integral involving
a holomorphic volume form on the fibers $Y_t$
at the $B$--side.
 
\smallskip
A fuller formulation of the mirror partnership
relation consists in the identification of two
{\it formal Frobenius manifolds (FM):} quantum cohomology
of any $X_s$ at the $A$--side, and Barannikov--Kontsevich's
FM on a formal {\it extended moduli space} at the $B$--side.

\smallskip

Most of the remaining part of this section will be devoted to the
description of the relevant Frobenius manifolds. 
\smallskip

However, 
the reader must be aware of more global aspects of this essentially local picture. In fact, moduli stack of complex variations of
$Y$ may have many cusps; they are acted upon by the 
Teichm\"uller group $\roman{Diff}\,(Y)/\roman{Diff}_0(Y)$. One can speculate
that mirror partnership is stable
with respect to such moduli cusp changes. Then the question
arises, what corresponds to them at the $A$--side.

\smallskip

Two partial answers were suggested. In [AsGM]
it was argued that different birational models
of some $X_s$ can produce canonically isomorphic cohomology
groups, in particular $H^2$, in which however K\"ahler cones will form
a non--trivial fan  (in the sense of toric geometry).
Maximal cones of this fan support K\"ahler cusps that
might correspond to different moduli cusps of the same
family at the $B$--side. In this picture, one does not
see what should correspond to the Teichm\"uller group.
M.~Kontsevich suggested in the framework of his conjectured
Lagrangian/complex duality that it must be the
autoequivalence group of the derived category of
coherent sheaves on $X_s$. Some evidence for this
was furnished by comparison of the stabilizing subgroups
of the cusps:  see [Hor], [SeTh].

\medskip

{\bf 1.3. Frobenius manifolds.} Let $M$ be an
analytic or formal supermanifold. A structure
of the Frobenius manifold on it is given
by a flat metric $g$ (symmetric non--degenerate
form on the tangent sheaf) and a function (potential)
$\Phi$ with the following property. Let $(x_a)$
be a local $g$--flat coordinate system, $\partial_a=
\partial /\partial x_a$, $\Phi_{abc}=\partial_a
\partial_b\partial_c\Phi$. Raise one index of $\Phi_{abc}$
using $g$ and define an $\Cal{O}_M$--bilinear 
multiplication $\circ$ on $\Cal{T}_M$ by
$\partial_a\circ\partial_b:=\sum_c \Phi_{ab}{}^c\partial_c.$
Then this multiplication must be associative
(it is obviously (super)commutative).
Additional structures that are present in the mirror picture
are the flat identity $e$ for $\circ$ and an Euler
vector field $E$ satisfying the conditions $\roman{Lie}_E(g)=Dg$
for some constant $D$ and
$\roman{Lie}_E(\circ )=\circ$. It expresses homogeneity
properties of $\Phi$: we have $E\Phi=(D+1)\Phi$ +
a polynomial in flat coordinates of degree $\le 2.$  

\smallskip

At the $A$--side, the relevant Frobenius manifold
is formal: $M$ is the formal completion
of the linear space $H^*(X_s,\bold{C})$, with its Poincar\'e
pairing as $g$ and the potential $\Phi$
constructed as formal series whose Taylor coefficients
are Gromov--Witten invariants of $X_s$.
At the $B$--side the relevant Frobenius manifold
can be conceived as a certain formal neighborhood
of the classical moduli space $T$ near the relevant cusp:
{\it extended moduli space} of the $B$--family.
Both formal spaces can be refined to germs
of analytic spaces.

\smallskip

Here are some details.

\medskip

{\bf 1.4. Quantum cohomology.} Potential $\Phi$ of the quantum
cohomology can be defined for any projective complex 
(or compact symplectic)
manifold $X$. After taking into account the relevant
homogeneity properties of the Gromov--Witten invariants 
it can be written as a formal series in linear coordinates
$(x^a)$ on $\oplus_{k\ne 2}H^k(X,\bold{C})$
and their exponentials on $H^2(X)$:
$$
\Phi (x)=\frac{1}{6}\,((\sum x_a\Delta_a)^3)
$$
$$
+
\sum_{\beta \ne 0}
e^{(\beta ,\sum_{|\Delta_{b}|=2} x_b\Delta_b)}
\sum_{n\ge 0, (a_i):\,|\Delta_{a_i}|\ne 2}
\langle \Delta_{a_n}\dots \Delta_{a_1} \rangle_{0,n,\beta}\,
\frac{x_{a_1}\dots x_{a_n}}{n!}.
\eqno(1.1)
$$
Here $(\Delta_a)$ is the basis of $H^*(X,\bold{C})$
dual to $(x_a)$, $\Delta_k \in H^{|\Delta_k|}(X)$,
the first term in the rhs of (1.1)
is the cubic self--intersection index, $\beta$ runs over
numerically effective 1--classes in $X$. Finally,
the Gromov--Witten invariant
$\langle \Delta_{a_n}\dots \Delta_{a_1} \rangle_{0,n,\beta}$
counts virtual number of stable maps of genus zero
$(C; x_1,\dots ,x_n; f:\,C\to X)$ such that
$f_*([C])=\beta$ and $f(x_i)\in D_{a_i}$ where $D_{a_i}$
is a cycle representing homology class dual to $\Delta_{a_i}.$
Physically, $\Delta_a$ are called
{\it the primary fields} of the respective
Conformal Field Theory, and the Gromov--Witten invariants are 
their {\it correlators}.

\smallskip

The {\it small quantum cohomology potential}
is obtained by restricting $\Phi (x)$ to $H^2$,
that is, putting $x_a=0$ for $|\Delta_a|\ne 2.$

\medskip

{\bf 1.5. Barannikov--Kontsevich's construction.}
On the $B$--side, the relevant formal Frobenius
potential is constructed on the completion at zero
of the cohomology space $H^*(Y,\wedge^*(\Cal{T}_Y))$ interpreted
as a formal moduli
space  $\Cal{M}_{\Cal{A}_{\infty}}$ of $A_{\infty }$--deformations of
$Y$. This construction was introduced in [Bar];
it refines the earlier proposal from [BK]. 
Unlike the case of quantum cohomology,
here it is essential to require $Y$ to be a (weak) Calabi--Yau
manifold. This condition will be used, in particular,
through a choice of the global holomorphic volume form
$\Omega$ on $Y$.

\smallskip

This geometric setup produces first of all an
algebraic object $(\Cal{A},\delta ,\Delta, \int )$,
{\it special differential Batalin--Vilkovyski algebra
(dBV)},
consisting of the following data which we will
describe in axiomatized form.

\smallskip

(i) {\it $\Cal{A}$ is a supercommutative $\bold{C}$--algebra.}

In the Calabi--Yau setup,
$\Cal{A}=\Gamma_{C^{\infty}}(Y,\wedge^*(\overline{T}_Y^*)
\otimes\wedge^*(T_Y))$.

\smallskip

(ii) {\it $\delta$ is an odd $\bold{C}$--derivation of $\Cal{A}$,
$\delta^2=0.$}

In our case, $\delta =\overline{\partial}$, the
operator defining the complex structure on $Y$ and
its tangent bundle, so that $\Cal{A}$ is the Dolbeault
resolution of the exterior algebra of the tangent bundle.

Therefore, the $\delta$--cohomology space
$H=H(\Cal{A},\delta )=\roman{Ker}\,\delta/\roman{Im}\,\delta$
in our case is identified with coherent cohomology
$H^*(Y,\wedge^*(\Cal{T}_Y)).$ Generally, we assume it to be
of finite dimension.
\smallskip

The space $H$ plays the central role, because
it will support the structure of the formal Frobenius manifolds.

\smallskip

We will denote by $K=\bold{C}[[x_a]]$  the ring
of formal functions on $H,$ $(x_a)$
being coordinates on $H$ dual to a basis
$(\Delta_a).$

\smallskip

(iii) {\it $\Delta$ is another odd differential, $\Delta^2=0$, which
is a differential operator of order two
with respect to the multiplication in $\Cal{A}.$}

More precisely, we assume that for any $a\in \Cal{A}$
the formula
$$
\partial_ab=(-1)^{\tilde{a}}\Delta (ab)- (-1)^{\tilde{a}}(\Delta a)b
-a\Delta b
$$
defines a derivation $\partial_a$. Moreover, we assume that 
$\delta\Delta +\Delta\delta =0$.

In our case, $\Delta$ is obtained from the
$\partial$--operator on the complexified
$C^{\infty}$ de Rham complex of $Y$ after the
identification of this complex with $\Cal{A}$
with the help of $\Omega$: 
$\Delta (a):= (\vdash \Omega )^{-1}\circ\partial\circ (a\vdash \Omega ).$
  
\smallskip

From the $\partial\overline{\partial}$--Lemma in K\"ahler geometry,
it follows that the two canonical embeddings
of differential spaces
$$
(\roman{Ker}\,\Delta ,\delta ) \to (\Cal{A},\delta ), \quad
(\roman{Ker}\,\delta ,\Delta ) \to (\Cal{A},\Delta )
\eqno(1.2)
$$
are quasi--isomorphisms, and moreover, homology
of all four differential spaces can be idemtified 
with $(\roman{Ker}\,\Delta\cap\roman{Ker}\,\delta )/\roman{Im}\,\delta\Delta .$

\smallskip

As a part of this package, one also obtains the following formality property:
the natural map $\roman{Ker}\,\Delta \to H(\Cal{A},\Delta )$
induces surjection of differential Lie algebras which is a quasi--isomorphism:
$$
(\roman{Ker}\,\Delta ,\,[\bullet ],\,\delta )\to (H(\Cal{A},\Delta ),0,0) .
$$

\smallskip

In the axiomatized situation, we impose these conditions 
as an additional axiom. This condition can be weakened: it suffices
to require only that  cohomology of differentials
$\delta +\Delta$ and $\delta$ have the same dimension.

\smallskip

\smallskip

(iv) {\it $\int :\Cal{A}\to \bold{C}$ is a linear functional
which must satisfy two integration by parts identities:}
$$
\int (\delta a)b= (-1)^{\widetilde{a}+1}\int a\delta b,\quad
\int (\Delta a)b= (-1)^{\widetilde{a}}\int a\Delta b .
\eqno(1.3)
$$

\smallskip

The integral is given by the formula
$$
\int a=\int_Y (a\vdash \Omega)\wedge\Omega
\eqno(1.4)
$$
where $\Omega$ means a holomorphic volume form on $Y$
whose period over the unique monodromy invariant
cycle at the chosen cusp is $(2\pi i)^d, d=\roman{dim}\,Y.$

\smallskip

(v) {\it Algebra grading $\Cal{A}=\oplus \Cal{A}^n,\, \bold{C}\in
\Cal{A}^0$.}

\smallskip

We assume that with respect to this grading, $\delta$ and $\Delta$
are of degree 1, and $\int$ has a definite
degree. (This is at variance with [Ma4], [Ma5], but agrees with
[Bar]).

\smallskip

Grading produces an Euler field on $H$, whereas the image
of $1\in \Cal{A}$ serves as flat identity.

\smallskip

In the Calabi--Yau setup, we can grade 
$\wedge^p\overline{T}^*_Y\otimes
\wedge^q T_Y$ by $q-p$.

\medskip

{\bf 1.5.1. Frobenius structure.} Having thus described the
formal properties of a Batalin--Vilkovyski algebra
$(\Cal{A},\delta ,\Delta, \int )$, we can now explain
the derivation of the Frobenius structure on $H$.

\smallskip

One starts with checking that the bilinear operation
$[a\bullet b]=\partial_ab$, together with
multiplication, endows $\Cal{A}$ by the structure
of {\it Gerstenhaber, or odd Poisson superalgebra}, in which the Lie
bracket is a parity changing operation, and all the
usual axioms are valid after inserting appropriate signs.

\smallskip

The basic ingredient of the construction from [Bar] is a certain exponential
map ${\Phi}^W $. In the Calabi-Yau setup it is an $A_{\infty}$--analog
$\Cal{M}_{A_{\infty}}
\rightarrow H^{*}(Y,\bold{C})[[\hbar ^{-1},\hbar ]][d] $
of the classical period map. Roughly speaking the map $ \Phi ^{W} $ is described
by the formula
$$
\Phi ^{W}(x_{a,}\hbar )=\left[ \exp \frac{1}{\hbar }\widetilde{\Gamma }\right]
$$
 where $ \widetilde{\Gamma }\in \Cal {A}\widehat{\otimes }K[[\hbar ]] $
is a $ W- $normalized generic solution to the Maurer-Cartan equation 
$ (\delta +\hbar \Delta )
\widetilde{\Gamma }+\frac{1}{2}[\widetilde{\Gamma }\bullet \widetilde{\Gamma }]=0 $
and $ \left[ a\right]  $ denotes the cohomology class with respect to the
differential $ \delta +\hbar \Delta  $. Here $ \delta  $ and $ \Delta  $
are assumed to be extended to $ \Cal {A}\widehat{\otimes }K[[\hbar ]] $
by linearity and $ \widetilde{\Gamma } $ is supposed to be $ W- $normalized
generic in the following sense: firstly, 
$ \left[ \exp \frac{1}{\hbar }\widetilde{\Gamma }\right] \in 1+L_{W} $,
where $ L_{W} $, $ \hbar ^{-1}L_{W}\subset L_{W} $ is semi-infinite subspace
associated with an increasing isotropic filtration on cohomology of $ \delta
+\Delta  $,
and, secondly, the map $ (\Phi ^{W}-1)\, \roman{mod}\, (\hbar ^{-1}L_{W}):
H\rightarrow L_{W}/\hbar ^{-1}L_{W} $
is linear and is an isomorphism. In the Calabi-Yau setting $ W $
is the monodromy weight filtration associated with the relevant cusp. Existence
of such solution $ \widetilde{\Gamma } $ for $ W $ satisfying certain
transversality condition can be proved by induction on the order of coefficients
of Taylor expansion. 

\smallskip

As a matter of fact, at this stage this construction exhibits certain common
features with the K.Saito's construction of FM structures on unfolding spaces
of singularities. It seems that if one chooses for  $W$ a certain special filtration
then the primitive form from the K.Saito theory can be identified with an analog
of $ \Phi ^{W}(x_{a,}\hbar ) $. The existence of a primitive
form in K.~Saito's theory is a nontrivial fact which follows
in general from the theory of mixed Hodge modules of M.~Saito.

\smallskip

Let us put now $ \Gamma =\widetilde{\Gamma }(x^{a},\hbar =0) $ and 
$ \delta _{\Gamma }:=id\otimes \delta +[\Gamma \bullet ] $.
The operator $ \delta _{\Gamma } $ is a homological differential acting on
$ \Cal {A}_{K}:=K\widehat{\otimes }\Cal {A} $. By continuity, one
can canonically identify $ H(\Cal {A}_{K},\delta _{\Gamma }) $ with $ K\otimes H $.
On the other hand, multiplication in $ \Cal {A}_{K} $ induces a multiplication
on $ H(\Cal {A}_{K},\delta _{\Gamma }) $. This is our $ \circ  $.
The map $ \Phi ^{W}(x_{a,}\hbar ) $ induces a pairing on the tangent sheaf
to $ H $: 
$$
\left\langle \partial _{a},\partial _{b}\right\rangle ^{W}:=
\int \partial _{a}\Phi ^{W}(x_{a,}\hbar )\partial _{b}\Phi ^{W}(x_{a,}-\hbar )
$$
The properties of the map $ \Phi ^{W} $ imply that this pairing is constant:
$ \left\langle \partial _{a},\partial _{b}\right\rangle ^{W}=g_{ab} $.
This is our flat metric.

\medskip

{\bf 1.5.2. Mirror identities
for complete intersections in projective spaces.} 
After these preparations, Barannikov's proof runs as
follows. Barannikov invokes the famous
Givental's result ([Giv2], [Giv5], [LiLY1]) establishing the mirror identity
on the level of "small quntum cohomology" (restriction to $H^2$)
replacing $A$--model, and classical moduli space replacing
$B$--model. This furnishes identification
of a part of Gromov--Witten invariants as coming from
the relevant Picard--Fuchs equations.
Now, Kontsevich--Manin's ``First reconstruction theorem''
from [KoM] shows that this part suffices
for the identification of the remaining invariants as soon as we know
that Associativity Equations (= Frobenius structure) hold.
In dimension 3 the latter supply no additional information, but
the larger dimension is, the more important
Associativity Equations become.

\medskip

{\bf 1.5.3. Extended moduli spaces.} The context of Mirror Symmetry
served to increase awareness of the importance of {\it extended moduli spaces}
in many other contexts of algebraic geometry.
Roughly speaking, any classical deformation problem 
is governed by a cohomology group $H^k$ classifying
infinitesimal extensions and the next cohomology group $H^{k+1}$
classifying obstructions. In the stable and unobstructed
case, $H^k$ is the tangent space to the base of versal deformation.
Extended moduli space in the unobstructed case
has total cohomology $H^*$ as tangent space. Barannikov--Kontsevich's
$B$--model is such an extended moduli space for Calabi--Yau manifolds.

\smallskip
See [KoS], [CiKa], [Mane] for a discussion of this matter
in general, and [Me3] for interesting constructions, related
to the Frobenius structure.

\medskip

{\bf 1.6. Other mirror isomorphisms.} There exist 
isomorphisms of auxiliary Frobenius manifolds
connecting certain unfolding spaces of singularities ($B$--model)
and moduli
spaces of curves with spin structure ($A$--model)
respectively, as was suggested by Witten [Wi2] and mathematically developed in [JaKV1], [JaKV2]. See also [Ma4]
about possible relations to the Calabi--Yau mirror picture,
developing the context in which the Mirror Symmetry
was first discussed in [Ge1], [Ge2].

\bigskip

\centerline{\bf 2. Lagrangian/complex duality and Mirror Symmetry}

\medskip

{\bf 2.1. Classical phase spaces.} Consider a $C^{\infty}$ symplectic
manifold $(X,\omega )$, endowed with a submersion
$p_X:\,X\to U$ whose fibers are Lagrangian tori,
and a Lagrangian section $0_X:\,U\to X.$ This is the classical
setup of action--angle variables in the theory of completely
integrable systems.

\smallskip

The form $\omega$ identifies the bundle of Lie algebras
of the tori $p^{-1}_X(u),\,u\in U,$ with the cotangent bundle
$T^*_U.$ Hence $T^*_U$ can be seen as fiberwise
universal cover of $X$, and we have a canonical isomorphism
$X=T^*_U/H$ where $H$ is a Lagrangian sublattice in $T^*_U$
with respect to the lift of $\omega$ which is the standard
symplectic form on the cotangent bundle. There exists
also a canonical flat symmetric connection
on $T^*_U$ for which $H$ is horizontal.

\smallskip

Put $H^t=\Cal{H}om\,(H,\bold{Z}).$ This local system is embedded as
a sublattice into $T_U$, and we can define {\it the mirror
partner} of $(p_X:\,X\to U,\omega ,0_X)$ as the toric fibration
$Y:=T_U/H^t$ endowed with the projection  to the same
base $p_Y:\,Y\to U$
and the zero section $0_Y$.

\medskip

{\bf 2.2. Complex structure on $Y$.} Passing from
$X$ to $Y$ we have lost the symplectic form. To compensate for this loss,
we have acquired a complex structure $J:\,T_Y\to T_Y$
which can be produced from $(p:\,X\to U,\omega ,0_X)$
in the following way. The flat connection on
$T_U$ obtained by the dualization from $T^*_U$
produces a natural splitting $T_Y=p_Y^*(T_U)\oplus p_Y^*(T_U).$
With respect to this splitting, $J$ acts as $(t_1,t_2)\mapsto (-t_2,t_1).$

\smallskip

Conversely, suppose that we have a complex manifold $Y$ endowed
with a fibration by real tori $Y\to U$ with zero section,
such that the operator of complex structure along the zero section
identifies $T_U$ with the bundle of Lie algebras of fibers.
Then we can consecutively construct the lattice $H^t\subset T_U$,
the dual fibration $X:=T^*_U/H$ and the symplectic form
on $X$ coming from the cotangent bundle.

\medskip

{\bf 2.3. Fourier--Mukai transform and further relationships
between Lagrangian and complex geometry.} Consider first
a pair of dual real tori $\bold{T}=H_{\bold{R}}/H$
and $\bold{T}^t=H_{\bold{R}}^t/H^t$ where $H$ is a free 
abelian group of finite rank, $H^t$ the dual group.
Denote by $\langle \,,\rangle$ the scalar product
$H^t\times H\to\bold{Z}$ and its real extensions.
Each point $x^t\in\bold{T}^t$ can be interpreted
as a local system of one dimensional complex vector spaces
with monodromy $\pi_1(\bold{T})=H\to S^1:\,h\mapsto
e^{2\pi i\langle x^t,h\rangle}$.
Hence $\bold{T}^t$ becomes the moduli space of
such systems on $\bold{T}$, and similarly
with roles of $\bold{T}$ and $\bold{T}^t$ reversed.

\smallskip

This can be conveniently expressed by introducing
{\it the Poincar\'e bundle $(\Cal{P},\nabla_{\Cal{P}})$}
on $\bold{T}\otimes\bold{T}^t$ which is rank one complex bundle with
connection. The connection is flat along both projections, but
has curvature $2\pi i\langle\partial^t,\partial\rangle$
on $(\partial^t,\partial ) \in H^t\times H.$

\smallskip

Using $(\Cal{P},\nabla_{\Cal{P}})$, we can extend the correspondence
between points of $\bold{T}$ and local systems on $\bold{T}^t$
in the following way. Call a skyscraper sheaf $\Cal{F}$ on $\bold{T}$
a sheaf consisting of a finite number of vector spaces $F_i$
supported by points $x_i.$ We can define a functorial map
$$
\Cal{F}\mapsto p_{\bold{T}^t*}( p_{\bold{T}}^*(\Cal{F})\otimes\Cal{P})
\eqno(2.1)
$$
whose image, if one takes in account the induced connection, 
is {\it a unitary local system on $\bold{T}^t$},
that is, a complex vector bundle with flat connection
and semisimple monodromy with eigenvalues in $S^1$.

\smallskip

Let now $X$ and $Y$
be mirror partners in the sense of 2.1--2.2. The construction above shows
first of all
that points $y$ of $Y$ bijectively correspond to pairs
consisting of a Lagrangian torus $L=p_X^{-1}(p_Y(y))$
and a unitary local system of rank one on it. 

\smallskip

Moreover, $X\times_U Y$ carries the relative Poincar\'e bundle
which we again will denote $(\Cal{P},\nabla_{\Cal{P}})$:
connection is extended in an obvious way in the horizontal directions.
An appropriate relative version of skyscraper sheaves is
played by pairs $(L,\Cal{L})$ consisting of a Lagrangian submanifold
of $X$ transversal to the tori and a unitary local system $\Cal{L}$
on $L$. The Fourier transform (2.1) of such a system is defined by
$$
(L,\Cal{L})\mapsto p_{Y*}(p_L^*\Cal{L}\otimes (i\times\roman{id})^*\Cal{P})
\eqno(2.2)
$$
where we denote by $i:\,L\to X$ the Lagrangian immersion,
and $p_Y:\,L\times_U Y\to Y,$ $i\times\roman{id}:\,L\times_UY\to Y,$
$p_L:\,L\times_UY\to L.$ The image of (2.2) also carries the induced
connection. We can calculate the $\overline{\partial}$--component
of it in the complex structure of $Y$ and find out that it is flat.
In other words, the rhs of (2.2) is canonically a holomorphic
vector bundle on $Y$.

\smallskip

{\bf 2.3.1. An example: mirror duality between complex
or $p$--adic abelian varieties.} In this subsection we
propose a definition of mirror duality for abelian varieties
which works uniformly well over arbitrary complete
normed fields $K$. We will represent
such a variety $\Cal{A}$ as a quotient (in the analytic category)
of an algebraic $K$--torus $T$ by a discrete subgroup
$B$ of maximal rank. Such a ``multiplicative uniformization''
goes back to Jacobi. The passage to the algebraic--geometric picture
is mediated by the classical or $p$--adic theta--functions
which are defined as analytic functions on $T$
with the usual automorphic properties with respects to
shifts by elements of $B$, see e.g. [Ma6] for details. The choice of multiplicative 
uniformization adequately models the choice
of a cusp in the moduli space of abelian varieties.

\smallskip

To be precise, algebraic torus $T$ with the character group $H$
over a field $K$ is the spectrum of the group ring of $H$.
The dual torus $T^t$, as above, has the character group $H^t$.

\smallskip

Consider now any diagram of the form
$$
(j,j^t):\ T(K)\leftarrow B\rightarrow T^t(K)
\eqno(2.3)
$$
where  $B$ is free abelian group of the same rank as $H$ and
$j$, resp $j^t$, are its embeddings as discrete subgroups
into $T(K),$ resp. $T^t(K)$.

\smallskip

We will say that pairs $(\Cal{A}:=T(K)/j(B), j^t)$ and
$(\Cal{B}:=T^t(K)/j^t(B), j)$ are mirror dual to each other. 
The quotient spaces $\Cal{A}$, $\Cal{B}$ not always have
the structure of abelian varieties, but this is
not important for the following.

\smallskip

In order to motivate this definition, we will show that for $K=\bold{C}$,
we can produce from (2.3) a pair of dual real toric fibrations
over a common base.

\smallskip

We have
the Lie group isomorphism $\bold{C}^*\to S^1\times\bold{R}:$
$z\mapsto (z/|z|, \roman{log}\,|z|)$. This induces
an isomorphism
$$
(\alpha ,\lambda ):\ T(\bold{C})\to \roman{Hom}\,(H,S^1)
\times \roman{Hom}\,(H, \bold{R}).
\eqno(2.4)
$$
Since $j(B)$ is discrete of maximal rank,
then $\lambda\circ j(B)$ is an additive lattice in the real space
$\roman{Hom}\,(H, \bold{R})$.
Thus (2.4) produces a real torus fibration of $T(\bold{C})$ 
over the base which
is as well a real torus of the same dimension:
$$
0\to \roman{Hom}\,(H,S^1)\to T(\bold{C})/j(B)\to
\roman{Hom}\,(H,\bold{R})/\lambda\circ j(B)\to 0 \,.
\eqno(2.5)
$$
Similarly, we have 
$$
0\to \roman{Hom}\,(H^t,S^1)\to T^t(\bold{C})/j^t(B)\to
\roman{Hom}\,(H^t,\bold{R})/\lambda^t\circ j^t(B)\to 0 
\eqno(2.6)
$$
where $\lambda^t$ is defined for $T^t$ in the same way
as $\lambda$ for $T$. Let us identify linear real spaces $H_{\bold{R}}$
with $H^t_{\bold{R}}$ in such a way that lattice points
$\lambda\circ j(b)$ and $\lambda^t\circ j^t(b)$ are identified
for all $b\in B$. Then (2.5) and (2.6) become
dual real torus fibrations over the common base. 

\smallskip

The relevant complex structures in our context come from covering tori.
They produce symplectic forms as was explained above.

\medskip

{\bf 2.4. Kontsevich's package.} We now return to the general
mirror dual  toric fibrations. With some stretch of
imagination, one can see the following pattern
in the picture described above: {\it Lagrangian
cycles with local systems on $X$, whose projection to $U$
have real dimension $k$, must correspond to
coherent sheaves on $Y$ with support of complex
dimension $k$.} 

\smallskip

Kontsevich in [Ko2] suggested a considerably more sophisticated
conjecture. Namely,
let $X$  be
a {\it compact} symplectic manifold with $c_1(X)=0$, and  $Y$ some compact complex Calabi--Yau manifold.

\smallskip

Then the relation of mirror partnership between $X$ and $Y$
consists in an equivalence between the Fukaya
triangulated category $D (Fuk_X)$ concocted out of
Lagrangian cycles with local systems on the one side, and 
(a subcategory of) $D^b(Coh_X)$
on the other side.

\smallskip

Briefly, to construct $D(Fuk_X)$ one proceeds in three
steps: first, one constructs an $A_{\infty}$--category
$Fuk_Y$, then one produces from it
another $A_{\infty}$--category of twisted complexes,
and finally, one passes to the homology category of the latter.

\smallskip

Objects $\Lambda =(L,\Cal{L},\lambda )$ of $Fuk_Y$ 
are Lagrangian submanifolds $L$ in $X$ with
unitary local systems $\Cal{L}$, endowed with a lifting $\lambda$
to the fiberwise universal cover of the Lagrangian
Grassmannian of $X$.

\smallskip

Morphism space between a pair of such objects 
admits a transparent description in the case when their
Lagrangian submanifolds $L_1,L_2$ intersect transversally.
In this case it is simply $\roman{Hom}\, (\Cal{L_1},\Cal{L}_2)$
in the category of sheaves {\t on $X$}. This space is $\bold{Z}$--graded
with the help of a construction using $\lambda$ and Maslov index. 

\smallskip

However, the
composition of morphisms is not at all the composition
of these morphisms of sheaves. In fact, a modification
of Floer's construction using summation over pseudoholomorphic 
parametrized discs in $X$ produces
a series of polylinear maps
$$
m_1:\,\roman{Hom}\, (\Lambda_1 ,\Lambda_2 )\to \roman{Hom}\, (\Lambda_1,\Lambda_2),
$$
$$ 
m_2:\,\roman{Hom}\, (\Lambda_1,\Lambda_2)\otimes \roman{Hom}\, (\Lambda_2,\Lambda_3)\to \roman{Hom}\, (\Lambda_1,\Lambda_2),
$$
and generally
$$
m_r:\,\roman{Hom}\, (\Lambda_1,\Lambda_2)\otimes\dots\otimes \roman{Hom}\, (\Lambda_{r-1},\Lambda_r)\to \roman{Hom}\, (\Lambda_1,\Lambda_r).
$$
If the respective sums converge, $m_1$ endows the graded
$\roman{Hom}$--spaces with the structure of a complex,
$m_2$ becomes the morphism of complexes, and higher
multiplications are interrelated by the
$A_{\infty}$--identities ensuring that the associativity
constraints for the composition of morphisms are
valid up to explicit homotopies.

\smallskip

For more detailed discussion, see [Ko2], [PoZ], [Fu2], 
and the literature quoted therein. In particular, the case
of elliptic curves is rather well understood thanks
to Polishchuk and Zaslow, and Fukaya started
treating abelian varieties and complex tori.

\smallskip

Both categories involved in the Kontsevich's conjecture
generally have non--trivial discrete symmetries, induced 
in the CY--context by monodromy
at the Lagrangian side and by derived correspondences at the complex side.
Thus some additional data have to be chosen in order to
pinpoint the expected functor. The awareness of symmetries 
led Kontsevich
to beautiful predictions about the correspondence
between monodromy actions and
automorphims of derived categories:
see [Hor], [SeTh], [Tho]. We mentioned these predictions
above, when we discussed the global properties of
the Frobenius partnership relations.

\smallskip

Kontsevich was vague about both the origin
of the equivalence functor and exact geometric
relation between $X$ and $Y$. 
One can interpret the picture described in 2.1--2.3
which emerged later
as a precise guess about the nature of several data
left implicit in Kontsevich's presentation:

\smallskip

{\it (i) The character of additional data to be chosen:
dual toric fibrations of $X$, $Y$ over a common base.} 

We will see below how this choice
at the complex side is related to the notion of cusp
of the relevant moduli space which we introduced
in the context of Frobenius mirror partnership.

\smallskip

{\it (ii) The structure of the restriction of the equivalence functor
acting on the simple objects: Fourier--Mukai transform corresponding
to the choice (i).}

\smallskip

\smallskip

With exception of the case of complex tori,
there is not much chance that $X$ or $Y$ would
admit a global fibration by real tori:
degenerate fibers are generally unavoidable, and their geometry
and influence on the global geometry of the mirror
picture are poorly understood. The case of $K3$--surfaces
offers some testing ground, because $K3$--surfaces
are hyperk\"ahler, and Lagrangian tori can be
transformed into a pencil of elliptic curves
by an appropriate rotation of the complex structure.

\smallskip

Recently M.~Kontsevich and A.~Todorov came up with a conjectural
limiting metric picture of the maximally degenerating family of
CY manifolds of dimension $d$ (private communication).
 Namely, fix a cohomology
class of K\"ahler forms and a moduli cusp.
Deform the complex structure by moving to the maximal
degeneration point, and the
Calabi--Yau metric in the chosen class by multiplying it by a real number
in such a way that the diameter of the space remains 1.

\smallskip

Todorov and Kontsevich expect that the limit $\Cal{X}$ in the
Hausdorff--Gromov sense of this family of metric spaces 
will be a real $d$--dimensional
manifold with a Riemannian metric which might have
singularities in codimension two. Moreover, the remnants of
the special real torus fibration consist in the
following additional data: affine structure
and a sublattice in the tangent bundle. In local affine coordinates, the metric
must be the second derivative of a convex function $H$,
and the volume form of the metric must be constant. 

\smallskip

Conjecturally, mirror dual family (endowed with appropriate
cusps) produces the same limiting metric space $\Cal{Y}=\Cal{X}$, but
with a different affine structure and sublattice
in the lattice bundle.

\medskip

{\bf 2.5. Mirror Symmetry between Calabi--Yau manifolds.}
Let now $X$, $Y$ be two $C^{\infty}$--manifolds each of which
is endowed by a symplectic form, real toric fibration over a common base,
and a complex structure, $(\omega_X,p_X,J_X)$
and $(\omega_Y,p_Y,J_Y)$ respectively. We will say
that they are related by Mirror Symmetry, if
$(X,p_X,\omega_X)$ is the mirror partner
of $(Y,p_Y,J_Y)$ and $(X,p_X,J_X)$ is the mirror partner
of $(Y,p_Y,\omega_Y)$ in the sense of Lagrangian/complex duality.
An example of this setup is described in 2.3.1.

\smallskip

The structures $J$ and $\omega$ at each side, of course,
can be related. The most rigid connection between
them is the presence of the Riemann metric $g$
producing the K\"ahler package $(J,\omega , g)$.
In the case of Calabi--Yau manifolds, the natural choice
is Yau's Ricci--flat metric $g$.

\smallskip

The program of [StYZ]  develops this setup, in particular,
supplying the topological and the metric characterization
of the basic toric fibrations. Namely, the cohomology
class of any toric fiber in $X$, resp. $Y$ must be the generator
of the cyclic group of invariant cycles in the middle
cohomology with respect to the local monodromy action
at the chosen cusp of moduli space. Moreover, non--degenerate
toric fibers (and other relevant Lagrangian submanifolds)
 must be not simply Lagrangian,
but {\it special} Lagrangian. This produces
a version  of Lagrangian geometry whose rigidity is
comparable to that of complex one, and makes it fit for
comparison with the complex picture: see [Gr1], [Gr2], [Ty1], [Ty2]
for many details.

\smallskip

It would be important to develop a version of Fukaya's
category in this rigid context
where the usual tools of homological algebra
might work better.

\medskip

{\bf 2.6. Motives in the looking glass.} One of the most
basic expressions of the Mirror Symmetry of the Calabi--Yau
manifolds is the existence of highly nontrivial isomorphisms
between their cohomology spaces: the relation of
mirror partnership between $X$ and $Y$ is expected to produce,
roughly speaking,
an isomorphism $H^*(X)\to H^*(Y).$

\smallskip

More precisely, any isomorphism between the quantum cohomology
of $X$ and Barannikov--Kontsevich formal Frobenius manifold of $Y$
produces an identification of their spaces of flat vector fields, that is
a mirror isomorphism of the cohomology spaces
$$
\mu_{X,Y} :\quad H^*(X,\bold{C})\to H^*(Y,\wedge^*(\Cal{T}_Y))\otimes V_Y^{-2},
\quad V_Y:=H^0(Y,\Omega_Y^{\roman{max}}).
\eqno(2.7)
$$
Near a cusp in the moduli space of $Y$, $V_Y$ can be trivialized
by the choice of a volume form $\Omega$ having period $(2\pi i)^{\roman{dim}\,Y}$
along the invariant cycle. Then (2.7) becomes
a ring isomorphism. Trace functionals and flat metrics on both sides are identified
via (1.4). Comparing Euler fields, one sees that
$H^{p,q}(X)$ is identified with $H^q(Y,\wedge^p(\Cal{T}_Y).$
In particular, $H^{1,1}(X)$ becomes $H^2(Y,\Cal{T}_Y)$,
and the induced integral structure on the latter space
(exponential coordinates near the cusp) are described in [De3].

\smallskip

Notice now that the Frobenius structure at the left hand
side of (2.7) is essentially motivic, in the sense that numerical
Gromov--Witten invariants of $X$ come from algebraic
correspondences between $X^n$ and $\overline{M}_{0,n}$, $n\ge 3.$
More generally, theory of Gromov--Witten invariants
can be conceived as a chapter of algebraic and/or non--commutative geometry
{\it over the category of motives,} replacing the
more common category of linear spaces. This geometry deals,
for example, with affine groups whose function rings are
Hopf algebras in the category of $Ind$--motives. P.~Deligne
developed basics of this geometry in [De1], [De2], in order
to clarify the notion of motivic fundamental group.
Further examples come from or are motivated by  physics:
besides Gromov--Witten invariants, one can mention
Nakajima's theory of Heisenberg algebras related to Chow schemes of surfaces,
and a recent paper [LosMa].

\smallskip

It makes sense to ask then, what can be the mirror reflection
of this motivic geometry. Since the mirror maps are highly
transcendental, developing the adequate language
presents an interesting challenge.
Starting with the category of motives
in the sense of [An]  
generated by Calabi--Yau manifolds,
we can try to extend it by adding mirror isomorphisms as
new motivated morphisms. In this context, Kontsevich's
correspondence between CY Teichm\"uller groups
and autoequivalences of derived categories might
have an analog, saying that the mirror isomorphisms
connect the motivic fundamental groups (see [De2]) and
motivic automorphism groups of CYs 
whose Lie algebras were studied in [LoLu]. For abelian varieties,
this phenomenon is stressed in [GolLO].

\medskip

{\it Acknowledgement.} I am grateful to S.~Barannikov,
M.~Kontsevich and Y.~Soibelman, who suggested
revisions and corrections to the first version
of this talk. Of course, I am fully responsible
for the final text.

\bigskip
\centerline{\bf Bibliography}

\medskip

[An] Y.~Andr\'e. {\it Pour une th\'eorie inconditionnelle
des motives.} Publ. Math. IHES, 83 (1996), 5--49.

\smallskip

[AP] D.~Arinkin, A.~Polishchuk. {\it Fukaya category
and Fourier transform.} Preprint   math.AG/9811023

\smallskip

[AsGM] P.~Aspinwall, B.~Greene, D.~Morrison. {\it Calabi--Yau spaces,
mirror manifolds and spacetime topology change in string theory.}
In: [Mirs2], 213--279.

\smallskip

[Bar] S.~Barannikov. {\it Extended moduli spaces and mirror
symmetry in dimensions $n>3.$} Preprint math.AG/9903124.

\smallskip

[BarK] S.~Barannikov, M.~Kontsevich. {\it Frobenius manifolds
and formality of Lie algebras of polyvector fields.} Int. Math. Res. 
Notices, 4 (1998), 201--215.

\smallskip

[Ba1] V.~Batyrev. {\it Dual polyhedra and the mirror symmetry
for Calabi--Yau hypersurfaces in toric varieties.} Journ. Alg. Geom.,
3 (1994), 493--535.

\smallskip

[Ba2] V.~Batyrev. {\it Variation of the mixed Hodge structure
of affine hypersurfaces in algebraic tori.} Duke Math. J.,
69 (1993), 349--409.

\smallskip

[Ba3] V.~Batyrev. {\it Quantum cohomology ring of toric manifolds.}
Ast\'erisque, 218 (1993), 9--34.

\smallskip

[BaBo1] V.~Batyrev, L.~Borisov.  {\it On Calabi--Yau complete
intersections in toric varieties.} In: Proc. of Int. Conf.
on Higher Dimensional Complex Varieties (Trento, June 1994),
ed. by M.~Andreatta, De Gruyter, 1996, 39--65.

\smallskip

[BaBo2] V.~Batyrev, L.~Borisov. {\it Dual cones and mirror symmetry
for generalized Calabi--Yau manifolds.} In: Mirror Symmetry II,
ed.~by S.~T.~Yau, 1996, 65--80.

\smallskip

[BaBo3] V.~Batyrev, L.~Borisov. {\it Mirror duality and
string--theoretic Hodge numbers.} Inv.~Math., 126:1 (1996), 183--203.

\smallskip

[BaS] V.~Batyrev, D.~van Straten. {\it Generalized hypergeometric
functions and rational curves on Calabi--Yau complete intersections
in toric varieties.} Comm. Math. Phys., 168 (1995), 493--533.

\smallskip

[BerCOV] M.~Bershadsky, S.~Cecotti, H.~Ooguri, C.~Vafa. {\it Kodaira--Spencer
theory of gravity and exact results for quantum string amplitudes.}
Comm. Math. Phys., 165 (1994), 311--427.

\smallskip

[Ber] A.~Bertram. {\it Another way to enumerate rational
curves with torus actions.} math.AG/9905159.

\smallskip

[BiCPP] G.~Bini, C.~de Concini, M.~Polito, C.~Procesi.
{\it On the work of Givental relative to Mirror Symmetry.}
  math. AG/9805097.

\smallskip

[Bor1] L.~Borisov. {\it On Betti numbers and Chern classes of
varieties with trivial odd cohomology groups.}  
alg-geom/9703023.

\smallskip

[Bor2] L.~Borisov. {\it Vertex algebras and mirror symmetry.}
math.AG/9809094.

\smallskip

[Bor3] L.~Borisov. {\it Introduction to the vertex algebra approach
to mirror symmetry.}   arXiv:math.AG/9912195.

\smallskip

[BorLi] L.~Borisov, A.~Libgober. {\it Elliptic genera of toric
varieties and applications to mirror symmetry.}  
math.AG/9904126 (to appear in Inv. Math.)

\smallskip

[CaOGP] P.~Candelas, X.~C.~de la Ossa, P.~S.~Green, L.~Parkes.
{\it A pait of Calabi--Yau manifolds as an exactly soluble
superconformal theory.} Nucl. Phys., B 359 (1991), 21--74.

\smallskip

[CK] E.~Cattani, A.~Kaplan. {\it Degenerating variations of Hodge structures.}
Ast\'e\-risque, 179--180 (1989), 67--96.

\smallskip

[Ce] S.~Cecotti. {\it $N=2$ Landau--Ginzburg vs. Calabi--Yau
$\sigma$--models: non--perturbative aspects.} Int. J. of Mod. Phys. A,
6:10 (1991), 1749--1813.

\smallskip

[CiKa] I.~Ciocan--Fontanine, M.~Kapranov. {\it Derived Quot schemes.}
Preprint  math.AG/9905174

\smallskip

[CoK] D.~Cox, S.~Katz. {\it Mirror symmetry and algebraic geometry.}
AMS, Providence RI, 1999.

\smallskip

[De1] P.~Deligne. {\it Le groupe fondamental de la droite
projective moins trois points.} In: Galois groups over $\bold{Q}$,
ed. By Y.~Ihara, K.~Ribet, J.~P.~Serre, Springer Verlag, 1999,
79--297.

\smallskip

[De2] P.~Deligne. {\it Cat\'egories tannakiennes.} In: The Grothendieck
Festschrift, vol. II, Birkh\"auser Boston, 1990, 111--195.

\smallskip

[De3] P.~Deligne. {\it Local behavior of Hodge structures at infinity}. 
In: Mirror Symmetry II,
ed.~by B.~Greene and S.~T.~Yau, AMS--International Press, 1996, 683--699.

\smallskip

[DeMi] P.~Deligne, J.~Milne. {\it Tannakian categories.}
In: Springer Lecture Notes in Math., 900 (1982), 101--228.

\smallskip

[Dol] I.~Dolgachev. {\it Mirror symmetry for lattice polarized
K3 surfaces.} J.~Math.~Sci. 81 (1996), 2599--2630.
  alg-geom/9502005

\smallskip 

[DoM] R.~Donagi, E.~Markman. {\it Cubics, integrable systems, and
Calabi--Yau threefolds.} In: Proc. of the Conf. in Alg. Geometry
dedicated to F. Hirzebruch, Israel Math. Conf. Proc., 9 (1996).

\smallskip

[Dor] Ch.~Doran. {\it Picard--Fuchs uniformization: modularity
of the mirror map and mirror--moonshine.}
  math.AG/9812162

\smallskip

[Du1] B.~Dubrovin. {\it Geometry of 2D topological field theories.}
In: Springer LNM, 1620 (1996), 120--348.

\smallskip

[Du2] B.~Dubrovin. {\it Geometry and analytic theory of Frobenius
manifolds.} Proc. ICM Berlin 1998, vol. II, 315--326.
  math.AG/9807034.

\smallskip

[Fu1] K.~Fukaya. {\it Morse homotopy, $A^{\infty}$--categories,
and Floer homologies.} In: Proc. of the 1993 GARC Workshop
on Geometry and Topology, ed. by H.~J.~Kim, Seoul Nat. Univ.

\smallskip

[Fu2] K.~Fukaya. {\it Mirror symmetry of abelian variety and
multi theta functions.} Preprint, 1998.

\smallskip

[Ge1] D.~Gepner. {\it On the spectrum of 2D conformal field theory.}
Nucl. Phys., B287 (1987), 111--126.

\smallskip

[Ge2] D.~Gepner. {\it Exactly solvable string compactifications
on manifolds of $SU(N)$ holonomy.} Phys. Lett. B199 (1987), 380--388.

\smallskip

[Ge3] D.~Gepner. {\it Fusion rings and geometry.} Comm. Math. Phys.,
141 (1991), 381--411.

\smallskip

[Giv1] A.~Givental. {\it Homological geometry I: Projective
hypersurfaces.} Selecta Math., new ser. 1:2 (1995), 325--345.

\smallskip

[Giv2] A.~Givental. {\it Equivariant Gromov--Witten invariants.}
Int. Math. Res. Notes, 13 (1996), 613--663.

\smallskip

[Giv3] A.~Givental. {\it Stationary phase integrals, quantum
Toda lattices, flag manifolds and the mirror conjecture.}
In: Topics in Singularity Theory, ed. by A.~Khovanski et al.,
AMS, Providence RI, 1997, 103--116.

\smallskip

[Giv4] A.~Givental. {\it Homological geometry and mirror symmetry.}
In: Proc. of the ICM, Z\"urich 1994, Birkh\"auser 1995, vol. 1,
472--480.

\smallskip

[Giv5] A.~Givental. {\it A mirror theorem for toric complete
intersections.} In: Topological Field Theory, Primitive Forms
and Related Topics, ed. by M.~Kashiwara et al., Progress in Math.,
vol. 60, Birkh\"auser, 1998, 141--175.   
alg--geom/9702016

\smallskip

[Giv6] A.~Givental. {\it Elliptic Gromov--Witten invariants and
the generalized mirror conjecture.}   math.AG/9803053.

\smallskip

[Giv7] A.~Givental. {\it The mirror formula for quintic threefolds.}
math.AG/9807070.

\smallskip

[GolLO] V.~Golyshev, V.~Lunts, D.~Orlov. {\it Mirror symmetry for
abelian varieties.}   math.AG/9812003

\smallskip

[Gre] B.~Green. {\it Constructing mirror manifolds.}
In: [MirS2], 29--69.

\smallskip

[Gr1] M.~Gross. {\it Special Lagrangian fibrations I: Topology.}
  alg-geom/9710006

\smallskip

[Gr2] M.~Gross. {\it Special Lagrangian fibrations II: Geometry.}
  math.AG/9809072

\smallskip

[GroW] M.~Gross, P.~M.~H.~Wilson. {\it Mirror symmetry via
3--tori for a class of Calabi--Yau threefolds.} Math. Ann.
309 (1997), 505--531.  
alg-geom/9608004.

\smallskip

[H] N.~Hitchin. {\it The moduli space of special Lagrangian
submanifolds.}   dg-ga/9711002.

\smallskip

[Hor] R.~P.~Horja. {\it Hypergeometric functions
and mirror symmetry in toric
varieties.}   arXiv:math.AG/9912109.

\smallskip

[HosLY1] S.~Hosono, B.~H.~Lian, S.~T.~Yau. {\it GKZ--generalized
hypergeometric systems in mirror symmetry of Calabi--Yau
hypersurfaces.}   alg--geom/9511001

\smallskip

[HosLY2] S.~Hosono, B.~H.~Lian, S.~T.~Yau. {\it Maximal
degeneracy points of GKZ systems.}  Journ. of AMS,
10:2 (1997), 427--443. alg--geom/9603014.



\smallskip

[JaKV1] T.~Jarvis, T.~Kimura, A.~Vaintrob. {\it The moduli
space of higher spin curves and integrable hierarchies.}
  math.AG/9905034

\smallskip

[JaKV2] T.~Jarvis, T.~Kimura, A.~Vaintrob. {\it Tensor
products of Frobenius manifolds and  moduli
spaces of higher spin curves.}
  math.AG/9911029



\smallskip

[Ko1] M.~Kontsevich. {\it $A_\infty$--algebras in mirror symmetry.}
Bonn MPI Arbeitstagung talk, 1993.

\smallskip

[Ko2] M.~Kontsevich. {\it Homological
algebra of Mirror Symmetry.} Proceedings of the ICM
(Z\"urich, 1994), vol. I, Birkh\"auser, 1995, 120--139.
  alg-geom/9411018.

\smallskip

[Ko3] M.~Kontsevich. {\it Mirror symmetry in dimension 3.} S\'eminaire Bourbaki, n${}^o$ 801, Juin 1995.

\smallskip

[Ko4] M.~Kontsevich. {\it Enumeration of rational curves via torus actions.}
In: The Moduli Space of Curves, ed. by
R\. Dijkgraaf, C\. Faber, G\. van der Geer, Progress in Math\.
vol\. 129, Birkh\"auser, 1995, 335--368.

\smallskip



[KoM] M.~Kontsevich, Yu.~Manin. {\it Gromov--Witten classes,
quantum cohomology, and enumerative geometry.} Comm. Math. Phys.,
164:3 (1994), 525--562.


[KoS] M.~Kontsevich, Y. Soibelman. {\it Deformation of
algebras over operads and Deligne's conjecture.}
math.QA/0001151

\smallskip

[KrS] M.~Kreuzer, H.~Skarke. {\it Complete classification of
reflexive polyhedra in four dimensions.}  
arXiv:hep-th/0002240.

\smallskip

[Ku] V.~S.~Kulikov. {\it Mixed Hodge structures and
singularities.} Cambridge Univ. Press, 1998.

\smallskip

[LYZ] N.~C.~Leung, Sh.-T.~Yau, E.~Zaslow. {\it From special
Lagrangian to Hermitian--Yang--Mills via Fourier--Mukai
transform.} Preprint math.DG/0005118.

\smallskip

[LiTY] B.~H.~Lian, A.~Todorov, S.~T.~Yau.
{\it Maximal unipotent monodromy for complete intersection
CY manifolds.} Preprint, 2000.

\smallskip

[LiLY1] B.~H.~Lian, K.~Liu, S.-T.~Yau. {\it Mirror principle I.}
Asian J. of Math., vol. 1, no. 4 (1997), 729--763.

\smallskip

[LiLY2] B.~H.~Lian, K.~Liu, S.-T.~Yau. {\it Mirror principle II.}
Asian J. of Math., vol. 3, no. 1 (1999), 109--146.

\smallskip

[Lib] A.~Libgober. {\it Chern classes and the periods of mirrors.}
  math.AG/9803119

\smallskip

[LibW] A.~Libgober, J.~Wood. {\it Uniqueness of the
complex structure on K\"ahler manifolds of certain homology
type.} J.~Diff.~Geom., 32 (1990), 139--154.

\smallskip

[LoLu] E.~Looienga, V.~Lunts. {\it A Lie algebra attached
to a projective variety.} Inv. Math., 129 (1997), 361--412.

\smallskip

[LosMa] A.~Losev, Yu.~Manin. {\it New moduli spaces of pointed curves and pencils
of flat connections.}
To be published in Fulton's Festschrift,
Michigan Journ. of Math., math.AG/0001003

\smallskip

[MalSV] F.~Malikov, V.~Schechtman, A.~Vaintrob. {\it Chiral
de Rham complex.} Comm. Math. Phys., 204:2 (1999), 439--473.
  math.AG/980341.

\smallskip

[Mane] M.~Manetti. {\it Extended deformation functors, I.}
  math.AG/9910071.

\smallskip

[Ma1] Yu.~Manin. {\it Problems on rational points and rational curves
on algebraic varieties.} In: Surveys of Diff. Geometry, vol. II,
ed. by C.~C.~Hsiung, S.~-T.Yau, Int. Press (1995), 214--245.

\smallskip

[Ma2] Yu.~Manin. {\it Generating functions in algebraic geometry
and sums over trees.} In: The Moduli Space of Curves, ed. by
R\. Dijkgraaf, C\. Faber, G\. van der Geer, Progress in Math\.
vol\. 129, Birkh\"auser, 1995, 401--418.

\smallskip

[Ma3] Yu.~Manin. {\it Sixth Painlev\'e equation, universal elliptic curve,
and mirror of $\bold{P}^2$.} AMS Transl. (2), vol. 186
(1998), 131--151.   alg--geom/9605010.

\smallskip

[Ma4] Yu.~Manin. {\it Three constructions of Frobenius manifolds:
a comparative study.} Asian J. Math., 3:1 (1999), 179--220
(Atiyah's Festschrift).   math.QA/9801006.

\smallskip

[Ma5] Yu.~Manin. {\it Frobenius manifolds, quantum cohomology, and moduli
spaces.} AMS Colloquium Publications, vol. 47, Providence, RI, 1999,
xiii+303 pp.

\smallskip

[Ma6] Yu.~Manin. {\it Quantized theta--functions.} In: Common
Trends in Mathematics and Quantum Field Theories (Kyoto, 1990), 
Progress of Theor. Phys. Supplement, 102 (1990), 219--228.

\smallskip

[McL] R.~C.~McLean. {\it Deformations of calibrated submanifolds.}, 
1996.

\smallskip

[Me1] S.~Merkulov. {\it Formality of canonical symplectic complexes
and Frobenius manifolds.} Int. Math. Res. Notes, 14 (1998), 727--733.

\smallskip

[Me2] S.~Merkulov. {\it  Strong homotopy algebras of a K\"ahler
manifold.} Preprint  math.AG/9809172

\smallskip

[Me3] S.~Merkulov. {\it Frobenius${}_{\infty}$ invariants
of homotopy Gerstenhaber algebras.} Preprint math.AG/0001007  

\smallskip

[MirS1] S.--T.~Yau, ed. {\it Essays on Mirror Manifolds.} International Press
Co., Hong Cong, 1992.

\smallskip

[MirS2] B.~Greene, S.~T.~Yau, eds. {\it Mirror Symmetry II.},
AMS--International Press, 1996.

\smallskip

[Mo1] D.~Morrison. {\it Mirror symmetry and rational curves
on quintic threefolds: a guide for mathematicians.}
J. AMS 6 (1993), 223--247.

\smallskip

[Mo2] D.~Morrison. {\it Compactifications of moduli spaces inspired
by mirror symmetry.} Ast\'erisque, vol. 218 (1993), 243--271.

\smallskip

[Mo3] D.~Morrison. {\it The geometry underlying mirror symmetry.}

\smallskip

[Pa] R.~Pandharipande. {\it Rational curves on hypersurfaces (after
A.~Givental)}. S\'eminaire Bourbaki, Exp. 848,
Ast\'erisque 252 (1998), 307--340. math.AG/9806133.

\smallskip

[Po1] A.~Polishchuk. {\it Massey and Fukaya products on elliptic curves.}
Preprint  math.AG/9803017.

\smallskip

[Po2] A.~Polishchuk. {\it Homological mirror symmetry 
with higher products.} Preprint  math.AG/9901025

\smallskip

[PoZ] A.~Polishchuk, E.~Zaslow. {\it Categorical mirror symmetry:
the elliptic curve.} Adv. Theor. Math. Phys.,
2 (1998), 443--470.   math.AG/980119.

\smallskip



[Schwa] J.~H.~Schwarz. {\it Lectures on Superstring and
M Theory Dualities.}   hep--th/9607201

\smallskip

[SeTh] P.~Seidel, R.~Thomas. {\it Braid group actions
on derived categories of coherent sheaves.}  
arXiv:math.AG/0001043.

\smallskip



[StYZ] A.~Strominger, S.--T.Yau, E.~Zaslow. {\it Mirror symmetry
is $T$--duality.} Nucl.~Phys. B 479 (1996), 243--259.



\smallskip

[Tho] R.~P.~Thomas. {\it Mirror symmetry and actions of braid groups on 
derived categories.}   arXiv:math.AG/0001044.

\smallskip

[Ty1] A.~Tyurin. {\it Special Lagrangian geometry and slightly deformed
algebraic geometry (SPLAG and SDAG)}.   math.AG/9806006

\smallskip

[Ty2] A.~Tyurin. {\it Geometric quantization and mirror symmetry.} Preprint
  \newline math.AG/9902027

\smallskip

[Va1] C.~Vafa. {\it Topological mirrors and quantum rings.}
In: Essays on Mirror Manifolds, ed. by Sh.--T. Yau,
International Press, Hong--Kong 1992, 96--119.

\smallskip

[Va2] C.~Vafa. {\it Extending mirror conjecture to Calabi--Yau with
bundles.}   hep--th/9804131.

\smallskip

[Va3] C.~Vafa. {\it Geometric Physics.} Proc. ICM
Berlin 1998, vol. I,  537--556.

\smallskip

[Voi1] C.~Voisin. {\it Sym\'etrie miroir.} Panoramas et
synth\`eses, 2 (1996), Soc. Math. de France.

\smallskip

[Voi2] C.~Voisin. {\it Variations of Hodge structure of
Calabi--Yau threefolds.} Quaderni della Scuola Norm. Sup. di Pisa,
1998.

\smallskip

[Wi1] E.~Witten. {\it Mirror manifolds and topological field theory.}
in: [MirS1], 120--159.

\smallskip

[Wi2] E.~Witten. {\it Algebraic geometry associated with
matrix models of two--dimensional gravity.}
In: Topological Models in Modern Mathematics (Stony Brook, NY, 1991),
Publish or Perish, Houston, TX (1993), 235--269.

\smallskip

[Za] E.~Zaslow. {\it Solitons and helices: the search for a 
Math-Physics bridge.} Comm. Math. Phys, 175 (1996), 337--375.

\smallskip

[Zh] I.~Zharkov. {\it Torus fibrations of Calabi--Yau hypersurfaces
in toric varieties and mirror symmetry.} Duke Math. J.,
101:2 (2000), 237--258.
  alg-geom/9806091

\enddocument